\documentclass{amsart}

\usepackage{amssymb}
\usepackage{amsmath}

\input xy
\xyoption{all}

\newtheorem{thm}{Theorem}[section]
\newtheorem{prop}[thm]{Proposition}
\newtheorem{cor}[thm]{Corollary}
\newtheorem{lm}[thm]{Lemma}

\newtheorem{ques}[thm]{Question}

\newtheorem{defn}[thm]{Definition}

\newcommand{\bpf}{\noindent {\em Proof.} }
\newcommand{\epf}{\qed \vspace{+10pt}}

\def\ni{\noindent}
\def\spec{\operatorname{Spec}}

\def\ps{\vspace{4pt}}

\def\sqr#1#2{{\vcenter{\vbox{\hrule height.#2pt
\hbox{\vrule width.#2pt height #1pt \kern #1pt \vrule 
width.#2pt}\hrule height.#2pt}}}}

\def\qed{{\hfil \break  \rightline{$\sqr74$}}}

\def\pic{{\rm Pic}}
\def\sym{{\rm Sym}}

\def\A{{\mathbb A}}
\def\C{{\mathbb C}}

\def\P{{\mathbb P}}

\def\Z{{\mathbb Z}}

\def\cX{{\mathcal X}}
\def\cH{{\mathcal H}}

\def\cO{{\mathcal O}}

\def\cO{{\mathcal O}}
\def\M{\overline{M}}

\begin{document}

\title[Sections over curves]{Rational connectivity and 
sections of families over curves}

\

\author[T. Graber]{Tom Graber}
\address{Mathematics Department, Harvard University,
1 Oxford st., Cambridge MA 02138, USA}
\email{graber@math.harvard.edu}
\thanks{T.G. was partially supported by an NSF postdoctoral fellowship.}

\author[J. Harris]{Joe Harris}
\address{Mathematics Department, Harvard University,
1 Oxford st., Cambridge MA 02138, USA}
\email{harris@math.harvard.edu}

\author[B. Mazur]{Barry Mazur}
\address{Mathematics Department, Harvard University,
1 Oxford st., Cambridge MA 02138, USA}
\email{mazur@math.harvard.edu}

\author[J. Starr]{Jason Starr}
\address{Mathematics Department, M.I.T., Cambridge MA 02139, USA}
\email{jstarr@math.mit.edu}
\thanks{J.S. was partially supported by NSF grant DMS-0201423.}

\date{\today}

\maketitle

\section{Introduction}

\subsection{Statement of results}

In \cite{GHS}, it is proved that a one-parameter family of rationally
connected varieties always has a rational section: explicitly, we have
the

\ps

\begin{thm}\label{oldmainth} Let $\pi : X \to B$ be a proper morphism 
of complex varieties, with $B$ a smooth curve. If the general fiber of
$\pi$ is \textbf{rationally connected}, then $\pi$ has a section
(recall a complex variety $V$ is called rationally connected if any
two points on $V$ lie on the image of some map $h:\P^1 \rightarrow
V$).
\end{thm}

The goal of this paper is to state and prove a converse to this
statement (c.f. Theorem~\ref{mainth} below).

\

We should first of all discuss what we mean by this, inasmuch as the
literal converse of Theorem~\ref{oldmainth} is clearly false. To this
end, let's focus on the question: under what circumstances does a
family $\pi : X \to B$ of varieties have the property that its
restriction to a general curve $C \subset B$ has a section?

\

This is certainly the case if the family $\pi : X \to B$ has a global
rational section. It is also the case by Theorem~\ref{oldmainth} if
the general fiber of $\pi$ is rationally connected, and by extension
it's the case if $X$ contains a subvariety $Z \subset X$ dominating
$B$ and whose fiber over a general point of $B$ is rationally
connected. (We can think of the case where the family $\pi : X \to B$
has a global rational section as a special case of this, a single
point being a rationally connected variety!) In this paper, we will
prove that in fact these are the only circumstances under which it may
occur.  To make this claim precise, we start by making the definition

\begin{defn}
Let $\pi : X \to B$ be an arbitrary morphism of complex varieties. By
a \textbf{pseudosection} of $\pi$ we will mean a subvariety $Z \subset
X$ such that the restriction $\pi|_Z : Z \to B$ is dominant with
rationally connected general fiber.
\end{defn}

Our main result is

\begin{thm}\label{mainth}
Let $B$ be any variety.  For every positive integer $d$ there exists a
bounded family $\cH_d$ of maps $h:C \rightarrow B$ from smooth
irreducible curves to $B$ such that for any proper morphism $\pi : X
\to B$ of relative dimension $d$ or less, if $h:C \to B$ is a map
parametrized by a very general point of $\cH_d$, the pullback
$$
\pi_C : X_C = X \times_B C \to C
$$
has a section if and only if $\pi$ has a pseudosection.

If $B$ is normal and quasi-projective, we can take $\cH_d$ to
be the family of smooth linear sections of $B$ under a sufficiently
positive projective embedding.
\end{thm}

What we are saying here, in other words, is that the class of
rationally connected varieties is the \textit{largest} class for which
the statement of Theorem~\ref{oldmainth} holds: if we have any other
family of varieties $\pi : X \to B$ satisfying the condition that
every one-parameter subfamily has a section, it does so by virtue of
the fact that $X$ contains a family of rationally connected varieties.

\

As a corollary of this, we will at the end of the paper settle a
question left hanging in [GHS]: whether or not the statement of
Theorem~\ref{oldmainth} holds for the larger class of
$\cO$-\textit{acyclic} varieties---that is, varieties $X$ with $H^i(X,
\cO) = 0$ for all $i > 0$. In fact, it was suggested by Serre in a
letter to Grothendieck ([GS], p.  152) that this might be the case
(though Serre immediately adds that it is ``sans doute trop
optimiste"). In Section~\ref{enriq}, we show this does not hold:
specifically, by applying Theorem~\ref{mainth} to the universal family
over a parameter space of Enriques surfaces with a particular
polarization, we will deduce the

\begin{cor}\label{enriqcor}
There exists a one-parameter family $X \to B$ of Enriques surfaces
with no rational section.
\end{cor}

\section{Stable maps and stable sections}

Our proof of Theorem~\ref{mainth} involves an induction on the
relative dimension of $f:X\rightarrow B$ where the base case
(i.e. fiber dimension zero) is proved by a version of the Lefschetz
hyperplane theorem.  In the course of the proof we will need to use
specializations of irreducible curves in $B$.  There are several
possible compactifications of the Chow variety of irreducible curves
in $B$, but the one we will use is the \emph{Kontsevich space of
stable maps}.  The reader who is unfamiliar with stable maps is
referred to the article \cite{FP}.  But the reader will not be far off
if he thinks of a stable map to $X$ as a map $h:C \rightarrow X$ for
which $C$ is a connected, projective curve which has at-worst-nodes as
singularities i.e. every point of $C$ is either smooth or is formally
isomorphic to the origin on the curve $\{(x,y)\in \A^2: xy = 0\}$.
There is an equivalence relation on stable maps which is the obvious
one, and the Kontsevich moduli space of stable maps is the
corresponding moduli space of equivalence classes of stable maps.

\

We let $\M_g(X,\beta)$ denote the Kontsevich space of stable maps
$h:C\rightarrow X$ such that $C$ has arithmetic genus $g$ and such
that the push-forward fundamental class $h_*[C]$ equals $\beta \in
H_2(X,\Z)$.  This is a proper Deligne-Mumford stack with projective
coarse moduli scheme (c.f. \cite{FP}).  We will often not need the
decorations, so let $\M(X) := \cup_{g,\beta} \M_g(X,\beta)$ denote the
space of all stable maps to $X$.

\begin{defn}If $h:C \to X$ and $h':C' \to X$ are
stable maps, we say that $h'$ is a \textbf{submap} of $h$ if there
is a factorization $h'=h\circ i$ for some $i:C' \to C$.
\end{defn}

Note in the above definition we do \emph{not} assume that $i$ is an
embedding.  The following definition makes precise what it means for
one stable map to be a submap of a specialization of another stable
map.

\begin{defn}

Given two irreducible families of stable maps, $D$ and $E$, we say
that $D$ \textbf{dominates} $E$ if every map parametrized by $E$ can be
realized as a submap of a specialization of maps parametrized by $D$.
Precisely, for any point $e \in
E$ parametrizing a stable map $h:C\rightarrow X$, we can find a family
of stable maps over $\spec (\C[[t]])$ such that $h$ is a submap of the
map over the special fiber, and the map over the generic fiber is
pulled back by a morphism from $\spec \left(\C((t)) \right)$ to $D$.
If $D$ and $E$ are reducible, we say $D$ dominates $E$ if every
component of $D$ dominates every component of $E$.
\end{defn}

We need a criterion for when a family of stable maps to a projective
variety is dominated by a family of embedded complete intersection
curves.  We will repeatedly make use of the following criterion.

\begin{lm}\label{smoothinglm}
Suppose $E$ is a family of stable maps to a normal, projective variety
$B$ such that for the general map $h: C \to B$ parametrized by $E$,
there is an open subset $U\subset B$ contained in the smooth locus of
$B$ and such that:
\begin{enumerate}
\item $h^{-1}(U)$ is dense in $C$,
\item $h^{-1}(U)$ contains all the nodes of $C$, and
\item the restricted map $h^{-1}(U) \rightarrow U$ is a closed embedding.

\end{enumerate}
Then $E$ is dominated by the family of one dimensional linear sections
of $B$ under a sufficiently positive projective embedding.

\end{lm}
  
\bpf 

Let $h:C\to B$ be a general element of our family.  We embed $B$ in a
projective space and find an integer $a$ such that the ideal of the
reduced image curve, $h(C)$ is generated by polynomials of degree $a$.
The sufficiently positive embedding required is the $a$th Veronese
reembedding.  Now we know that we can realize $h(C)$ as an
intersection of hyperplanes.  If we choose $\dim(B)-1$ generic
hyperplane sections of $B$ which contain $h(C)$, their intersection
will be a generically reduced curve $C'$ which contains $h(C)$ as a
subcurve.  Choosing the hyperplanes so that they generate the ideal of
$h(C)$ at each of the images of the nodes of $C$, and after possibly
shrinking $U$ (so that it still satisfies our hypotheses), we can
arrange that $h^{-1}(U) \rightarrow U \cap C'$ is an isomorphism.

\

Now we choose any one parameter family of smooth complete
intersections in $B$ whose flat limit is $C'$, i.e. a morphism from $\spec
(\C[[t]])$ to the Hilbert scheme of complete intersections in $B$
whose general fiber maps to a smooth complete intersection and whose
special fiber maps to $C'$.  We think of the general fiber as a stable
map and perform stable reduction to the corresponding map $\spec
(\C[[t]]) \rightarrow \M(B)$.  Denote the special fiber of the stable
reduction by $\overline{h}:\overline{C} \rightarrow B$.
Since $U\cap C'$ is already
at-worst-nodal and stable (since it is an embedding), there is an open
subset $V\subset \overline{C}$ such that $\overline{h}:V \rightarrow U
\cap C'$ is an isomorphism.  In other words, we have a factorization
$i:h^{-1}(U) \rightarrow \overline{C}$ of $h:h^{-1}(U) \rightarrow
B$.  Since every point of $C$ in $C-h^{-1}(U)$ is smooth, we can apply
the valuative criterion of properness to extend this factorization to
a morphism $i:C \rightarrow \overline{C}$.  So $h:C \rightarrow B$ is
a submap of $\overline{h}:\overline{C} \rightarrow B$, which shows
that $E$ is dominated by the family of smooth curves in $B$ which are
complete intersections of $d-1$ hyperplanes. 

\epf

We remark that the lemma above is not the most general result, but to
prove a stronger version would lead us too far astray.  We leave
it to the interested reader to prove that in the above lemma it
suffices to assume that for the general map $h:C\rightarrow B$
parametrized by $E$, the preimage of the smooth locus,
$h^{-1}(B_{\text{smooth}})$, is a dense open set which contains every
node of $C$.

\

Although our main concern is to understand sections of a map $\pi : X
\to B$ over smooth curves in $B$, the specialization methods we use
force us to consider the more general notion of stable maps.
Similarly, we need to replace the notion of section over a curve with
an object which specializes as we specialize the base curve.  To
explain this notion, we first recall a standard construction.

\

Given a morphism of projective schemes $\pi : X \to B$, there is an
induced map on Kontsevich spaces, $\M (\pi) : \M(X) \to \M(B)$
(technically we must restrict to those stable maps with $g>1$ or with
$\beta \neq 0$, but this will always be the case for us).  This map
takes the moduli point $f : C\to X$ to the point $\pi \circ f : \tilde
C \to B$, where $\tilde C$ is the stabilization of $C$ relative to the
morphism $\pi \circ f$.  It is obtained by contracting those
components of $C$ of genus 0 which are contracted by $\pi \circ f$ and
meet the rest of $C$ in fewer than 3 points.

\

\begin{defn} Given a morphism $\pi: X \to B$, and
a stable map $h:C \to B$, we define a \textbf{stable section} of $\pi$
over $[h]$ to be a stable map $\tilde{h}:\tilde{C} \to X$ such that
$\M(\pi)([\tilde h])=[h].$
\end{defn}

Notice that for a stable section $\tilde h$, the class
$\pi_*\left(\tilde{h}\right)_*\left[\tilde{C}\right]$ is just $h_*[C]$
and $g(C') = g(C)$.  Since $h$ is a stable map, either $g(C) > 1$ or
$h_*[C]$ is nonzero, i.e. the map $\M(\pi)$ really is defined in a
neighborhood of $[\tilde{h}]$.

\

If $C$ is a smooth curve in $B$, then a stable section of $\pi$ over
$C$ is simply a section of $\pi$ over $C$ with some trees of
$\pi$-contracted rational curves in $X$ attached.  The notion is more
interesting when $C$ has nodes.  In this case, the existence of a
stable section over $C$ does not guarantee the existence of any
sections over $C$ whatsoever.  For example, a typical stable section
$\tilde{h}$ over a curve $C=C_1 \cup C_2$, where $t_1\in C_1$ is glued
to $t_2\in C_2$, would consist of sections $\tilde{h}_1$ and
$\tilde{h_2}$ of $\pi$ over $C_1$ and $C_2$ separately together with a
tree of $\pi$-contracted rational curves in $X$ joining
$\tilde{h}_1(t_1)$ to $\tilde{h}_2(t_2)$ (and some trees of
$\pi$-contracted rational curves attached elsewhere).  Such rational
curves are exactly the sort which are contracted under the
stabilization process associated with $\M(\pi)$.

\

The point of this definition is that given a family of curves in the
base $B$ specializing to some stable map $h:C\rightarrow B$ (possibly
reducible), and given an honest section over the generic curve in this
family, then we cannot conclude the existence of a section of $\pi$
over $h$, but we do conclude the existence of a stable section of
$\pi$ over $h$.  In other words, the existence of a stable section is
preserved under specialization.  This follows immediately from the
properness of the irreducible components of $\M(X)$.  Another
elementary fact is that if $h'$ is a submap of $h$, then the existence
of a stable section over $h$ implies the existence of a stable section
over $h'$.  Taken together, we get the following

\begin{lm}\label{specializationlm} 
If $D$ and $E$ are families of stable maps to $B$, with $D$ dominating
$E$ and if a general map parametrized by $D$ admits a stable section,
then so does every map parametrized by $E$.
\end{lm}

As an application, we can strengthen the easy direction of our
main theorem.

\begin{prop}  If $\pi : X \to B$ is a morphism of projective
varieties with $B$ smooth and if $\pi$ admits a pseudosection, then
for any smooth curve $C$, and for any morphism $h: C \to B$, the
pullback family $X_C \to C$ admits a section.
\end{prop}

\bpf 

Since any map from a smooth curve factors through the normalization of
its image, it suffices to prove this statement for maps birational
onto their image.  Let $Z$ be a pseudosection of $\pi$.  We already
know by Theorem~\ref{oldmainth} that the proposition is true for any
smooth curve such that the general fiber of $Z$ over the curve is
rationally connected.  In particular it holds for a generic complete
intersection curve in $B$.  By Lemma~\ref{smoothinglm}, our map $h$
can be realized as a submap of a limit of such curves.  Then
Lemma~\ref{specializationlm} implies that $\pi$ admits a stable
section over $[h]$.  Since $C$ is smooth, this implies that $X_C$
admits a section over $C$.  \epf

\section{Proof of main theorem}

In our proof of Theorem~\ref{mainth}, we begin by assuming that $B$ is
normal and that both $X$ and $B$ are projective.  After handling this
``special'' case, we give the (easy) argument which reduces the
general case to the special case.

\

We will prove the theorem by induction on the relative dimension of
$X$ over $B$.  We start with the case of relative dimension zero.

\begin{prop}\label{dim0}
Let $B\subset \P^n$ be a normal variety and $\pi : X \to B$ a
generically finite proper morphism.  Then $\pi$ admits a rational
section if and only if $\pi$ admits a section when restricted to a
general one dimensional linear section of $B$.
\end{prop}

\bpf

Let $B_0 \subset B$ be an open subset such that the restriction of
$\pi$ to $X_0 = \pi^{-1}(B_0)$ is finite and \'etale of some degree
$k$.  Choose a one dimensional linear section of $B$ which is a smooth
curve $C$, such that the natural map $\pi_1(C\cap B_0) \to \pi_1(B_0)$
is surjective.  Such a $C$ exists by \cite{gm}.  In fact, any smooth
linear section which meets $B \setminus B_0$ suitably transversally
will do.  Define $C_0=C \cap B_0$, and choose any point $c_0 \in C_0$.
Choosing an ordering of the points in $\pi^{-1}(c_0)$ gives us a
natural monodromy representation $\rho : \pi_1(B_0,c_0) \to S_k$ where
$S_k$ is the symmetric group on $k$ letters.  The monodromy
representation associated to the \'etale cover $\pi^{-1}(C_0) \to C_0$
is just the composition $\rho \circ i_*$ where $i : C_0 \to B_0$ is
the inclusion.  The statement that $\pi$ admits a section over $C$ is
equivalent to this \'etale cover of $C_0$ admitting a section, which
is equivalent to asking for the image of $\rho \circ i_*$ to be
contained in the stabilizer of an element.  Since $i_*$ is surjective,
this is equivalent to asking for the image of $\rho$ to be contained
in the stabilizer of an element.  This in turn is equivalent to the
existence of a section of $X_0$ over $B_0$, i.e. a rational section of
$\pi$.  

\epf

In handling the case of positive relative dimension, one of the main
ingredients needed is the following bend-and-break lemma for sections.

\begin{lm}\label{bendbreaklm}
Let $\pi : X \to C$ be a morphism with $C$ a smooth curve.  Let $p\in
X$ be an arbitrary point.  If there is a positive dimensional family
of sections of $\pi$ passing through $p$, then there is a rational
curve in $X$ passing through $p$ which is contracted by $\pi$.
\end{lm}

\bpf

Let $q=\pi(p)$.  Suppose we have a one parameter family of sections
passing through $p$.  This gives us a rational map $f:B\times C \to X$
over $C$ whose restriction to a general fiber $\{b\}\times C$ is a
section passing through $p$.  Suppose, by way of contradiction, that
$f$ is a regular morphism in a neighborhood of $B \times \{q\}$.
Since $f$ contracts $B \times \{q\}$, by the rigidity lemma it also
contracts $B\times \{c\}$ for all $c\in C$, i.e. our family is
constant which contradicts that it is positive dimensional.  Hence,
$f$ is not regular near $B\times \{q\}$.  So there is at least one
point of indeterminacy in $B \times \{q\}$.  We may form the minimal
blow-up of $B\times C$ necessary to resolve the indeterminacy locus of
$f$.  The exceptional divisor of this blow-up is a tree of rational
curves which intersects the proper transform of $B \times \{q\}$ and
which is mapped to a tree of $\pi$-contracted rational curves in $X$.
Therefore some rational curve in the exceptional divisor maps to a
$\pi$-contracted rational curve which meets $p$.  

\epf

We will apply this lemma in two ways.  The first application is to get
a uniform bound on the dimensions of spaces of sections.  First we
need a definition.

\begin{defn} If $\pi : X \to B$ is a morphism, we define the
  \textbf{rational curve locus}, $V(\pi)$, to be the union of all
$\pi$-contracted rational curves in $X$.  
\end{defn}

Let us pause to describe what sort of object $V(\pi)$ is.  Fixing an
ample divisor $H$ on $X$, then for each integer $d$ there is a finite
type Chow variety (or Hilbert scheme, or Kontsevich space of stable
maps) parametrizing $\pi$-contracted rational curves whose $H$-degree
is $d$.  Over each such Chow variety, there is a universal family of
$\pi$-contracted rational curves along with a map to $X$.  And
$V(\pi)$ is simply the union over all $d$ of the image of this map
from the universal family to $X$.  We think of $V(\pi)$ as just a set,
but we constantly use the fact that $V(\pi)$ is a set which is a
countable union of subvarieties of $X$.

\

Here is our uniform bound on the dimensions of spaces of sections.

\begin{lm}\label{dimensionbound}
If $\pi : X \to C$ is a morphism of relative dimension $d$
and if $\Sigma$ is a family of sections of $\pi$
such that a general section parametrized by $\Sigma$
is not contained in $V(\pi)$, then $\dim (\Sigma) \leq d$.
\end{lm}

\bpf

Choose a very general point $c$ of $C$. Let $ev_c : \Sigma \to X_c$ be
the map which evaluates a section at $c$.  Our hypotheses ensure that
$\dim (X_c) \leq d$ and that $ev_c(\Sigma)\not\subset V(\pi)$.
Lemma~\ref{bendbreaklm} then implies that $ev_c$ is generically finite
onto its image, yielding the desired bound.  

\epf

\begin{defn}
Let $L$ denote the family of smooth linear sections of $B$ under some
projective embedding.  We define the family of \textbf{triangles},
$T(L)\subset \M(B)$ to be the family of all stable maps which look
like a triangle with sides in $L$. That is, we consider the family of
morphisms from a curve of the form $C=C_1 \cup C_2 \cup C_3$ such that
every pair of components of the domain curve meet in a single point
(and these intersection points are pairwise distinct), and such that
each $f|_{C_i}$ $i=1,2,3$ is an embedding onto a curve parametrized by
$L$.
\end{defn}

The parameter space $T(L)$ is irreducible, and also the total space of
the universal family of curves over $T(L)$ is irreducible.  Note also
that a general stable map parametrized by $T(L)$ satisfies the
hypotheses of Lemma~\ref{smoothinglm}.

\

Our second application of Lemma~\ref{bendbreaklm} is in the
proof of the following lemma, which is the main step in the proof of 
Theorem~\ref{mainth}.

\begin{lm}\label{trianglelm}
Let $\pi:X \to B$ be a projective morphism to a normal variety, and
let $p \in X$ be any point such that $p \notin V(\pi)$.  Suppose also
that $p$ is not contained in the closure of the image of any rational
section of $\pi$ (notice that if $q$ is a smooth point of $B$, then
every rational section passing through $p$ is actually regular at
$q=\pi(p)$ since there are no $\pi$-contracted rational curves meeting
$p$).  Then a very general triangle passing through $q=\pi(p)$ admits
no stable sections passing through $p$.
\end{lm}

\bpf

It suffices to exhibit a single triangle with this property.  We will
show directly that this holds for a very general triangle with a
vertex at $q=\pi(p)$.  Choose a subfamily $H\subset L$ of curves
passing through $q$, such that for general $b \in B$ a finite (but
positive) number of members of $H$ pass through $b$.  We construct a
subset $\Omega \subset X$ which is a countable union of subvarieties
of $X$ in the following way.  For every finite type family of sections
of $\pi$ over curves in $H$ and which take the value $p$ at $q$, we
have a map from the base of this family to $\M(X)$.  Form the closure
of the image of this map, and define $\tilde\Omega$ to be the
countable union of all such closed subvarieties of $\M(X)$ arising from
the countably many Chow varieties of sections as above.  Notice that
$\tilde\Omega$ is not necessarily quasi-compact, but it is a closed
subvariety of $\M(X)$ (which is also not quasi-compact).

\

We can restrict the universal curve of $\M(X)$ over $\tilde\Omega$,
and there is a map from the total space of this universal curve to
$X$.  We define $\Omega$ to be the image of this map, so $\Omega$ is a
countable union of closed subvarieties of $X$.  Let $\tilde\Omega_0$
be any irreducible component of $\tilde\Omega$ and let $\Omega_0
\subset X$ be the (closed) image of $\tilde\Omega_0$.

\

Consider the restricted morphism $\M(\pi):\tilde\Omega_0 \rightarrow
\M(B)$.  The general point of $\Omega_0$ parametrizes a section over a
member of $H$, so the image of $\Omega_0$ under $\M(\pi)$ is contained
in the closure $\overline{H}$ of $H$.  And by Lemma~\ref{bendbreaklm},
the morphism 
$\M(\pi):\Omega_0 \rightarrow \overline{H}$ is generically finite.  
Therefore the map from the universal curve over $\Omega_0$ to the
universal curve over $\overline{H}$ is generically finite.  By
construction, the evaluation morphism from the universal curve over
$\overline{H}$ to $B$ is generically finite.  So finally we conclude
the restricted morphism $\pi:\Omega_0 \rightarrow B$ is generically
finite, i.e. for a general point $b\in B$ there are only finitely many
preimages of $b$ in $\Omega_0$.  Moreover, each of these finitely many
preimages lies on an honest section over a curve in $H$ passing
through $b$.  So in fact $\pi:\Omega_0 \rightarrow B$ is
unramified over $b$.
Since $\Omega$ is the union of
countably many sets $\Omega_0$, we conclude that for a very general
point $b\in B$ there are only countably many preimages of $b$ in
$\Omega$.  Moreover, if we choose $b$ very general, then every
irreducible component $\Omega_0$ whose image $\pi(\Omega_0)$ contains
$b$ actually surjects to $B$.  Since also $\pi:\Omega_0 \rightarrow B$
is unramified over $b$, we conclude that $\pi:\Omega_0 \rightarrow B$
is \'etale over $b$.  Thus for $b\in B$ a very general point, 
for \emph{every} irreducible component $\tilde\Omega_0$
of $\tilde\Omega$, the restricted map $\pi:\Omega_0\rightarrow B$ is
\'etale over $b$ -- possibly for the trivial reason that $b$ is not
contained in $\pi(\Omega_0)$.

\

Next we observe that any rational section $\rho:B \to \Omega$ includes
$p$ (i.e. the closure of the image of $\rho$ contains $B$).  First
observe that $\rho$ factors through one of the subsets
$\Omega_0\subset \Omega$.  Now for a general point $b$ in $B$,
$\pi:\Omega_0 \rightarrow B$ is unramified over $b$.  Therefore
$\rho:B \rightarrow \Omega_0$ is actually regular in a neighborhood of
$b$.  And the image $\rho(b)$ lies on some honest section $\tilde h:C
\rightarrow X$ over a curve $C$ in $H$ which contains $q$ and $b$ and
such that $\tilde h(q)=p$.  Since $\pi:\Omega_0 \rightarrow B$ is
unramified at $\rho(b)=\tilde h(b)$, we have that $\tilde h:C
\rightarrow X$ and $\rho|_C: C \rightarrow X$ are equal as rational
maps.  By the valuative criterion of properness, we conclude that
$\tilde h:C \rightarrow X$ factors through the closure of the image of
$\rho$, in particular $p=\tilde h(q)$ lies on the closure of the image
of $\rho$.  One of our hypotheses is that $p$ does not lie on the
closure of the image of any rational section.  Therefore we conclude
that there is no rational section of $\pi$ whose image is contained in
$\Omega$.

\

By the last paragraph, for each irreducible piece $\Omega_0$ of
$\Omega$, we have that $\Omega_0 \rightarrow B$ admits no rational
section.  Moreover, over a very general point $b\in B$, every one of
the maps $\pi:\Omega_0 \rightarrow B$ is \'etale.  So by the proof of
Lemma~\ref{dim0}, for a very general curve $C_3$ in $L$ passing
through a very general point $b\in B$ (and not necessarily passing
through $q$), for each irreducible piece $\Omega_0$ of $\Omega$, there
is no section of $\pi:\pi^{-1}(C_3) \rightarrow C_3$ whose image lies
in $\Omega_0$.  So there is no section of $\pi:\pi^{-1}(C_3)
\rightarrow C_3$ whose image lies in $\Omega$.

\

Since $p$ is not in $V(\pi)$, no section through $p$ can be contained
in $V(\pi)$.  It follows that for a very general point $b$ in $B$,
$\Omega \cap V(\pi)\cap\pi^{-1}(b) = \emptyset$.  Choose a very
general curve $C_3$ in $L$ as above, and choose two very general
points $r$ and $s$ on $C_3$.  Then $\Omega\cap \pi^{-1}(r)$ is a
countable set, disjoint from $V(\pi)$, and every point in this set
lies on an honest section over a curve in $H$ passing through $r$.
Given any point in this set, there are at most countably many sections
of $\pi$ over $C_3$ which take this value at $r$.  Hence, there is a
countable collection of sections of $\pi$ over $C_3$ whose value at
$r$ is contained in $\Omega$.  Any such section is not contained in
$\Omega$, and thus meets $\Omega$ in at most countably many points
with countably many images in $C_3$.  Choosing $s$ not to lie in any
of these countably many countable sets, we conclude that for any
section $\tilde{h}_3$ of $\pi$ over $C_3$ such that $\tilde{h}_3(r)$
is contained in $\Omega$, we have that $\tilde{h}_3(s)$ is \emph{not}
contained in $\Omega$.

\

Now we take our triangle to be $C=C_1 \cup C_2 \cup C_3$ where $C_1$
and $C_2$ are members of $H$ which join $q$ to $r$ and $q$ to $s$
respectively.  By way of contradiction, suppose there is a stable
section $\tilde{h}$ of $\pi$ over $C$ whose image contains $p$.  As we
have discussed, such a stable section consists of honest sections
$\tilde{h}_1$, $\tilde{h}_2$ and $\tilde{h}_3$ over $C_1$, $C_2$, and
$C_3$ respectively, along with some trees of $\pi$-contracted rational
curves attached which connect $\tilde{h}_1(r)$ to $\tilde{h}_3(r)$,
which connect $\tilde{h}_2(s)$ to $\tilde{h}_3(s)$ and which connect
$\tilde{h}_1(q)$ and $\tilde{h}_2(q)$ to $p$, if these points don't
already coincide.  By the definition of $\Omega$, the images
$\tilde{h}_1(C_1)$ and $\tilde{h}_2(C_2)$ are necessarily contained in
$\Omega$.  Since $r$ and $s$ were chosen to be very general, there are
no $\pi$-contracted rational curves over $r$ or $s$ which meet
$\Omega$, in particular, there is no tree of $\pi$-rational curves
which meets either $\tilde{h}_1(r)$ or $\tilde{h}_2(s)$.  So we must
have $\tilde{h}_1(r) = \tilde{h}_3(r)$ and $\tilde{h}_2(s) =
\tilde{h}_3(s)$.  Since $p$ is not contained in $V(\pi)$, also there
is no tree of $\pi$-contracted rational curves which meets $p$.
Therefore $\tilde{h}_1(q) = \tilde{h}_2(q) = p$.  The upshot is that,
after pruning any extraneous trees of $\pi$-contracted rational
curves, we have that $\tilde{h}$ is an honest section of $\pi$ over
the reducible curve $C$.

\
 
But now we have our contradiction: we have seen that for any section
$\tilde{h}_3$ of $\pi$ over $C_3$ such that $\tilde{h}_3(r)$ is
contained in $\Omega$, then $\tilde{h}_3(s)$ is not contained in
$\Omega$.  On the other hand we have by the last paragraph that
$\tilde{h}_3(r) = \tilde{h}_1(r)$ is contained in $\Omega$ and also
$\tilde{h}_3(s) = \tilde{h}_2(s)$ is contained in $\Omega$.  Therefore
we conclude there is no stable section $\tilde{h}$ of $\pi$ over $C$.

\epf

Of course this lemma doesn't tell us much in case the fibers of $\pi$
are uniruled.  Thanks to a construction of Koll\'ar-Miyaoka-Mori and
using Theorem~\ref{oldmainth}, we can always reduce to the case that
the fibers of $\pi$ are non-uniruled.

\begin{defn} Given a morphism $\pi : X \to B$, 
the \textbf{relative mrc fibration}
$$\xymatrix{
X \ar@{-->}[r]^\phi \ar[d]_\pi& W\ar[ld]^{\pi'}\\
B & }
$$
is a dominant rational map $\phi : X \to W$ of varieties over $B$ such
that a general fiber of $\phi$ is rationally connected and a general
fiber of $\pi'$ is not uniruled.
\end{defn}

The existence of the relative mrc fibration is established in
~\cite[theorem IV.5.9]{K}, although the equivalence of the definition
given there with the one above requires Theorem~\ref{oldmainth}.

\

Before applying Lemma~\ref{trianglelm} to our main theorem, we note a
corollary which is interesting in its own right.

\begin{cor} 
If $B$ is a normal, quasi-projective variety, and $\cH$ is any family
of curves in $B$ which dominates $T(L)$ and whose general member is
smooth, then for any projective morphism $\pi : X \to B$, the
following two conditions are equivalent:

\begin{enumerate}
\item A general point of $X$ lies in a pseudosection.
\item For a general curve $C$ parametrized by $\cH$, a general
point of $X_C = \pi^{-1}(C)$ lies on a section of $\pi:X_C \to C$.
\end{enumerate}
Note that the existence of such a family $\cH$ is ensured by
Lemma~\ref{smoothinglm}.
\end{cor}

\bpf

One direction follows from a stronger version of
Theorem~\ref{oldmainth} which states that if $\pi : Z \to C$ has
rationally connected general fiber there is a section of $\pi$ through
a general point of $Z$ (cf. \cite[IV.6.10]{K}).  The other direction
is more interesting.  If the general fiber of $\pi$ is not uniruled
then this follows from Lemma~\ref{trianglelm}, since in this case $X$
being generically covered by pseudosections is equivalent to it being
generically covered by rational sections.  Therefore, suppose that the
general fiber of $\pi$ is uniruled.  Consider the relative mrc
fibration of $\pi:X \rightarrow B$,
$$\xymatrix{
X \ar@{-->}[r]^\phi \ar[d]_\pi& W\ar[ld]^{\pi'}\\
B & }
$$
The two conditions we are trying to prove equivalent are both
preserved under blowing up of $X$ (to see this for the second
condition, restrict to the general case that $C$ is smooth).  Thus we
are free to blow up $X$ in order to resolve the indeterminacy locus of
$\phi$.  So we may assume that $\phi$ is actually a regular morphism.
Suppose that there is a section through a general point of $X_C$ over
a general curve $C$. Composing with $\phi$ gives us a section of $W_C$
through a general point of $W_C$ over a general curve $C$.  Since the
fibers of $\pi'$ are not uniruled, we conclude by the last paragraph
that a general point of $W$ is contained in a rational section of
$\pi'$.  Taking the preimage of this rational section under $\phi$
gives us a pseudosection of $\pi$ and passing through a general point.

\epf

Our proof of Theorem~\ref{mainth} proceeds similarly.  First we will
prove the result for maps whose general fiber is not uniruled, and
then we will handle the general case by appealing to the relative mrc
fibration.  In addition we will use an induction on the relative
dimension of $X$ over $B$.  We have already considered the case of
fiber dimension zero in Proposition~\ref{dim0}, thus suppose that $d >
0$.  By way of induction, assume that we have already constructed a
family $\cH_{d-1}$ of smooth curves in $B$ which cover $B$ and such
that for any morphism $\pi : X \to B$ of relative dimension less than
$d$, $\pi$ admits a section when restricted over a very general curve
in $\cH_{d-1}$ if and only if $\pi$ admits a pseudosection.  We
construct $\cH_d$ as follows.  First we construct a family of
reducible nodal curves by letting $\tilde \cH_d$ be the family of maps
$f:C\to X$ of the form $C= C_0 \cup C_1 \cup \cdots \cup C_{d+1}$,
where $[f_{C_0}]$ is a member of $\cH_{d-1}$ and the other $C_i$,
$i=1,\dots, d+1$ are triangles which are pairwise disjoint and which
each meet $C_0$ in a single node which is embedded in the smooth locus
of $B$.  Now take $\cH_d$ to be any family of smooth curves which
dominates $\tilde \cH_d$.  By Lemma~\ref{smoothinglm}, we can take
$\cH_d$ to be the family of linear sections of $B$ under a
sufficiently positive projective embedding.

\

We need to check that $\cH_d$ satisfies the desired property.  Namely,
suppose $\pi: X \to B$ is a projective morphism of relative dimension
less than or equal to $d$ which does not admit a pseudosection.  Then
we need to show that over a very general member of $\cH_d$, $\pi$ does
not admit a section.  By Lemma~\ref{specializationlm}, it suffices to
check that over a very general member of $\tilde\cH_d$, $\pi$ does not
admit a stable section.

\

First we will consider the case where the general fiber of $\pi$ is
not uniruled.  In order to later handle the uniruled case, it will be
useful for us to prove a statement that seems stronger than necessary.
As above, we let $V(\pi)$ be the rational curve locus which is the
union of all $\pi$-contracted rational curves in $X$, which is a
subset of $X$ which is a countable union of subvarieties.  We let
$Y(\pi)$ be the union of $V(\pi)$ and all rational sections of $\pi$.
This is also a countable union of subvarieties of $X$.  Note that
Lemma~\ref{trianglelm} says exactly that for any point $p$ in
$X\setminus Y(\pi)$, a very general triangle through $\pi(p)$ admits
no stable sections containing $p$.

\

\begin{lm}\label{technicallm}
Any stable section of $\pi$ over a very general member of
$\tilde\cH_d$ has values over $C_0$ contained in $Y(\pi)$.
\end{lm}

\

Before proving the lemma, we remark that (given the inductive
hypothesis) it immediately implies our theorem in the case where the
fibers are not uniruled.  If $\pi : X \to B$ is a morphism whose
general fiber is not uniruled and which does not admit a pseudosection
then $Y(\pi) = V(\pi)$ is a countable union of {\em proper}
subvarieties of $X$, i.e. it is a countable union of subvarieties
$Y(\pi)_0$ of $X$ such that the fiber dimension of $Y(\pi)_0
\rightarrow B$ is strictly less than $d$. By the induction assumption,
for a very general $C_0$ in $\cH_{d-1}$, there can be no honest
section of $\pi$ contained in any of the subvarieties $Y(\pi)_0$.
Thus, by the lemma, there can be no stable section of $\pi$ over $C =
C_0 \cup C_1 \cup \dots \cup C_{d+1}$.

\bpf

We now prove the lemma.  We imagine assembling our very general member
of $\tilde \cH_d$ one component at a time.  Pick a very general $C_0
\in \cH_{d-1}$ and let $\Sigma$ denote the parameter space of all
sections of $\pi$ over $C_0$ which are not contained in $Y(\pi)$.
This is the complement of a countable union of subvarieties 
in a countable union of subvarieties of $\M(X)$.   Denote the irreducible
components of $\Sigma$ by $\Sigma^\alpha$, and
by Lemma~\ref{dimensionbound} we conclude that each 
$\Sigma^\alpha$ has dimension less than or equal to $d$.

\

Our strategy now is simple.  The condition that a section over $C_0$
extends to a section over $C_0 \cup C_i$ should impose a condition by
Lemma~\ref{trianglelm}, and so after imposing $d+1$ conditions there
should be no sections left.  To prove this, we consider the chain
$$\Sigma_{d+1} \subset \Sigma_d \subset \cdots \subset \Sigma$$ 
where $\Sigma_i$ is defined to be the subset of $\Sigma$ parametrizing
sections of $\pi$ over $C_0$ which are not contained in $Y(\pi)$ and
which can extend over $C_0 \cup C_1 \cup \cdots \cup C_i$.  That is,
if we let $q_1, \ldots, q_{d+1}$ be the very general points at which
we attach the triangles, $\Sigma_i$ parametrizes those sections of
$\pi$ over $C_0$ whose value at $q_j$ agrees with the value of some
stable section of $\pi$ over $C_j$ for all $j\leq i$.  We will prove
by induction on $i$ that $\dim(\Sigma_i) \leq d-i$ for each $i$, in
particular $\Sigma_{d+1}$ is empty.

\

We have already seen that every component of $\Sigma=\Sigma_0$ has
dimension at most $d$, so this establishes the base case $i=0$.  By
way of induction, assume that every component of $\Sigma_k$ has
dimension at most $d-k$.  Now we want to show the result for $k+1$.
For any one of the countably many irreducible components
$\Sigma_k^\alpha$ of $\Sigma_k$, for a general point $q_{k+1}$ some
section $\tilde{h}^\alpha_0$ parametrized by $\Sigma_k^\alpha$ maps
$q_{k+1}$ to a point not contained in $Y(\pi)$, i.e. the point
$p^\alpha = \tilde{h}_0(q_{k+1})$ is not in $Y(\pi)$.  So if we choose
a very general point $q_{k+1}$ we can arrange that for every
irreducible component $\Sigma^\alpha_k$ of $\Sigma_k$, there is a
section $\tilde{h}^\alpha_0$ in $\Sigma^\alpha_k$ such that
$p^\alpha=\tilde{h}^\alpha_0(q_{k+1})$ is not contained in $Y(\pi)$.

\

Now for each $\alpha$, for a general triangle $C_{k+1}$ through
$q_{k+1}$, we conclude by Lemma~\ref{trianglelm} that there is no
stable section over $C_{k+1}$ which passes through $p^\alpha$.  So if
we choose a very general triangle $C_{k+1}$, we can arrange that for
every $\alpha$, there is no stable section over $C_{k+1}$ which passes
through any of the points $p^\alpha$.  So none of the sections
$\tilde{h}^\alpha_0$ extends to a stable section over $C\cup C_{k+1}$.
So for each $\alpha$, $\Sigma_{k+1}\cap \Sigma_k^\alpha$ is a proper
closed subvariety and thus has dimension strictly less than
$\text{dim}(\Sigma_k^\alpha) \leq d-k$.  Since we have
$$\Sigma_{k+1} = \cup_\alpha \left( \Sigma_{k+1} \cap \Sigma_k^\alpha
\right)$$ 
we conclude that every irreducible component of $\Sigma_{k+1}$ has
dimension at most $d-k-1$.  So the claim is proved by induction on
$k$.  In particular, we conclude that $\Sigma_{d+1}=\emptyset$,
i.e. over $C_0$ every section of $\pi$ which can be extended to stable
sections over $C$ is contained in $Y(\pi)$.  \epf

\

As discussed above, Lemma~\ref{technicallm} proves the induction step
in case the fibers of $\pi$ are not uniruled.  So to finish the
inductive proof of Theorem~\ref{mainth}, we are left to consider the
case where the fibers of $\pi$ are uniruled.  We argue by
contradiction.

\

By way of contradiction, assume that we have a morphism $\pi : X \to
B$ with no pseudosection, but which admits a section when restricted
to every element of $\cH_d$.  Let $\phi : X \to W$ be the relative mrc
fibration.  We may resolve the indeterminacy locus of $\phi$ by
blowing up:
$$\xymatrix{
\tilde X \ar[d]_\psi \ar[dr]^{\tilde \phi}& \\
X \ar@{-->}[r]^{\phi \quad} \ar[d]_\pi& W\ar[dl]^{\pi'}\\
B &}
$$
Let $E\subset X$ denote the fundamental locus of the morphism $\psi$
(i.e. the image under $\psi$ of the exceptional divisor of $\psi$).

\

Notice that the relative dimension of $\pi|_E:E\rightarrow B$ is
strictly less than $d$.  Suppose that $\pi|_E:E\rightarrow B$ admits a
section when restricted over a very general curve $C$ in $\cH_d$.  By
Lemma~\ref{specializationlm}, we conclude that $\pi|_E:E\rightarrow B$
admits a stable section when restricted over a stable map in
$\tilde\cH_d$.  In particular, since every curve in $\cH_{d-1}$ occurs
as the $C_0$-submap of a stable map in $\tilde\cH_d$, we conclude that
$E\rightarrow B$ admits a section when restricted over a very general
curve $C_0$ in $\cH_{d-1}$.  By the induction hypothesis, this implies
that there is a pseudosection of $\pi|_E:E \rightarrow B$.  But, in
particular, this implies there is a pseudosection of $\pi:X
\rightarrow B$ which contradicts our assumption.  So we conclude that
for a very general curve $C$ in $\cH_d$, $\pi_E:E \rightarrow B$
admits no section when restricted over $C$.  On the other hand, our
assumption is that $\pi:X \rightarrow B$ does admit a section over
$C$, i.e. there exists a section over $C$ which is not contained in
$E$.  This is the same as a rational section of $\pi\circ\psi:\tilde X
\rightarrow B$ over $C$.  Since $C$ is smooth, by the valuative
criterion of properness this rational section of $\pi\circ \psi$
extends to a regular section of $\pi\circ \psi$ over $C$.

\

Thus we find that $\pi \circ \psi$ admits a section over a very
general curve
$C$ in $\cH_d$.  Now if $\pi \circ \psi$ admits a pseudosection, so
does $\pi$ by simply taking the image of the pseudosection under
$\psi$.  Therefore we conclude that $\pi\circ\psi: \tilde X
\rightarrow B$ admits no pseudosection, but it does admit a section
when restricted over a very general curve in $\cH_d$.  Therefore, as
far as deriving a contradiction is concerned, we can replace
$X$ by $\tilde X$.  So from now on we  
assume that $\phi:X \rightarrow W$ is a regular morphism.

\

Let $W' \subset W$ denote the locus of points over which the fiber of
$\phi$ is not rationally connected.  Any rational section of $\pi'$
not contained in $W'$ gives rise to a pseudosection of $\pi$, which
doesn't exist by hypothesis.  Therefore all rational sections of
$\pi'$ are contained in $W'$.  Applying Lemma~\ref{technicallm}, we
find that over a general member $C$ of $\tilde\cH_d$, any stable
section of $\pi'$ maps $C_0$ into the subset $W'\cup V(\pi')$.  Thus,
any stable section of $\pi$ over $C$ has maps $C_0$ into the subset
$\phi^{-1}\left(W'\cup V(\pi')\right)$.

\

On the other hand, $\phi^{-1}\left(W'\cup V(\pi')\right)$ is a
countable union of proper subvarieties of $X$, each of which has
relative dimension at most $d-1$ over $B$.  So by the induction
hypothesis, every section of $\pi:X \rightarrow B$ over a very general
curve in $\cH_{d-1}$ has image which is not contained in
$\phi^{-1}\left(W'\cup V(\pi')\right)$.  Of course $C_0$ is a very
general curve in $\cH_{d-1}$, and so admits \emph{no} sections in this
locus.  So we conclude that over a very general member of
$\tilde\cH_{d}$, $\pi$ admits no stable section.  By
Lemma~\ref{specializationlm}, $\pi$ admits no section over a very
general member of $\cH_d$, and this is a contradiction of our
assumptions.

\

So, in the case that the fibers of $\pi$ are uniruled, we have
established the induction step by contradiction.  We had previously
established the induction step in the case of non uniruled fibers.
This finishes the proof of the induction step, so the main theorem is
proved by induction.  

\epf

\section{The general case}\label{general}

In the last section we proved the main theorem in case $B$ is normal
and quasi-projective and $\pi:X \rightarrow B$ is projective.  In this
section we will show how to reduce the general case to this case.  We
proceed by induction on the dimension.

\

Suppose that $B$ is a finite type algebraic variety.  Then by Chow's
lemma we can find a projective, birational morphism $B_1 \rightarrow
B$ such that $B_1$ is quasi-projective.  Now by Noether normalization,
the normalization $B_2 \rightarrow B_1$ of $B_1$ is a finite morphism.
Thus $f:B_2 \rightarrow B$ is a projective, birational morphism such
that $B_2$ is quasi-projective and normal.  Let $\cH_d$ be the family
of curves $C$ in $B_2$ constructed in the last section.  The
restriction of $f$ to a general curve in this family -- let's call
this restriction $h:C \rightarrow B_2$ -- is a nonconstant morphism,
i.e. it is a stable map.  Therefore, replacing $\cH_d$ by a Zariski
dense open subset, we may consider $\cH_d$ to be a family of stable
maps $h:C\rightarrow B$ with smooth domain.  The claim is that
Theorem~\ref{mainth} holds for $B$ and $\cH_d$.  We will prove this by
induction, but before proceeding to the induction argument we
introduce a little more notation.

\

Suppose that $\pi:X \rightarrow B$ is a proper morphism of relative
dimension at most $d$ which admits no pseudosection.  We need to prove
that for a very general map $h:C\rightarrow B$ in $\cH_d$, $\pi$
admits no section over $h$.  The base change $\pi_2:X\times_B B_2
\rightarrow B_2$ is a proper morphism of relative dimension at most
$d$ which admits no pseudosection, since the image under $\pi_1:X
\times_B B_2 \rightarrow X$ of a pseudosection of $\pi_2$ is a
pseudosection of $\pi$.  By again applying Chow's lemma, we can find a
projective, birational morphism $\phi:X_2 \rightarrow X$ such that
$\pi_2\circ\phi:X_2 \rightarrow B_2$ is projective.  Any pseudosection
of $\pi_2\circ \phi$ maps under $\phi$ to a pseudosection of $\pi_2$.
Therefore $\pi_2\circ\phi$ admits no pseudosection.

\

Now $\pi_2\circ\phi:X_2 \rightarrow B_2$ satisfies the hypotheses of
the last section.  By the proof of the main theorem in that section,
for a very general curve $C$ in $\cH_d$, $\pi_2\circ\phi$ admits no
section over $C$.  Let $Z\subset X$ denote the fundamental locus of
the birational, projective morphism $X_2 \rightarrow X$, i.e. the
locus over which this morphism is not an isomorphism.

\

If $d=0$, we are essentially done.  The locus $Z\subset X$ is a proper
subvariety, and since $\pi$ is generically finite, also $\pi(Z)\subset
B$ is a proper subvariety.  If we choose a very general map
$h:C\rightarrow B$ in $\cH_0$, then the image $h(C)$ does not lie in
$\pi(Z)$.  But then any section of $\pi:X \rightarrow B$ over $h$
determines a rational section of $\pi_2\circ\phi: X_2 \rightarrow B_2$
over $C$.  Since $C$ is smooth, by the valuative criterion of
properness this rational section extends to a regular section.  This
contradicts the result of the last section.  So we conclude that for a
very general map $h:C\rightarrow B$ in $\cH_0$, there is no section of
$\pi:X \rightarrow B$ over this map.

\

Now we proceed by induction.  We have established the base case $d=0$,
so we suppose that $d>1$.  By way of induction, we suppose the theorem
has been proved for $d-1$.  Consider $\pi|_Z:Z \rightarrow B$.  This
morphism has fiber dimension at most $d-1$.  By our induction
assumption, we conclude that $\pi|_Z:Z \rightarrow B$ has no section
when restricted over a very general map $h_0:C_0 \rightarrow B$ in
$\cH_{d-1}$.  By Lemma~\ref{specializationlm}, we conclude that
$\pi|_Z:Z\rightarrow B$ has no section when restricted over a very
general map $h:C \rightarrow B$ in $\cH_d$ (since $\cH_d$ dominates
$\cH_{d-1}$).  So if we choose a very general map $h:C\rightarrow B$
in $\cH_d$, then for any section $\tilde{h}:C \rightarrow X$ of $\pi$
over $h$, we have that $\tilde{h}(C)$ is not contained in $Z$.  So the
regular section $\tilde{h}$ determines a rational section of
$\pi_2:X_2 \rightarrow B_2$ over $C$.  Since $C$ is smooth, by the
valuative criterion of properness this rational section extends to a
regular section.  This contradicts the result of the last section.  So
we conclude that for a very general map $h:C\rightarrow B$ in $\cH_d$,
there is no section of $\pi:X \rightarrow B$ over $h$.

\section{Application: families of Enriques surfaces}\label{enriq}

In this section we'll show how to apply Theorem~\ref{mainth} to a
family of Enriques surfaces to deduce Corollary~\ref{enriqcor}, that
is, to find a one-parameter family of Enriques surfaces without a
section.

\subsection{A family of quartic Enriques surfaces}

The family we'll be starting with is the universal family over a
parameter space for quartic Enriques surfaces: that is, a family of
polarized Enriques surfaces $S$ with a polarization $M \in \pic(S)$ of
self-intersection 4 that includes a general such surface.  Now, for
the purposes of applying Theorem~\ref{mainth} and deducing
Corollary~\ref{enriqcor}, we can just write down the family as in
Definition~\ref{princenriq} below -- we don't need to know that it is
actually the generic quartic Enriques surface, and the reader who
doesn't particularly care can jump directly to
Definition~\ref{princenriq} -- but since we're going to be working
closely with the family it seems worthwhile to take a few paragraphs
and establish its bonafides.

\

To begin with, since Enriques surfaces $S$ have fundamental group
$\pi_1(S) \cong \Z/2\Z$ and have as universal covering space a K3
surface, a quartic Enriques surface $S$ is the quotient of an octic K3
surface $T$ -- that is, a K3 surface $T$ with a polarization $L \in
\pic(T)$ of self-intersection $c_1(L)^2 = 8$ -- by an involution
$\tau$ of $T$ preserving $L$.  For a generic octic K3 $(T,L)$, the
linear system of sections of $L$ is base-point-free and defines an
embedding of $T$ into $\P^5$, and the image surface is the
intersection of three quadric hypersurfaces in $\P^5$ with defining
equations $Q_1$, $Q_2$, and $Q_3$.

\

Next, since $\tau^*L \cong L$, the action of $\tau$ can be lifted to
an action on $H^0(T,L)$, and hence to an involution of $\P^5$ carrying
$T$ to itself. Moreover, since by Riemann-Roch
$$
h^0(S,M) \; = \; \frac{c_1(M)^2}{2} + \chi(\cO_S) \; = \; \frac{4}{2}
+ 1 \; = \; 3
$$
the action of $\tau$ on $H^0(T,L)$ must have eigenvalues $1$ and $-1$,
each with multiplicity 3. We thus have a canonical direct-sum
decomposition
$$
H^0(T,L) \; = \; \Gamma \oplus \Psi
$$
with $\dim \Gamma = \dim \Psi = 3$.

\

Applying the same principle, we see that the action of $\tau$ on
$H^0(T,L^2)$ has eigenvalue 1 with multiplicity
$$
h^0(S,M^2) \; = \; \frac{c_1(M^2)^2}{2} + \chi(\cO_S) \; = \; 
\frac{16}{2} + 1 \; = \; 9
$$
and correspondingly eigenvalue $-1$ with multiplicity $h^0(T,L^2)-9 =
18-9 = 9$. On the other hand, given that $H^0(T,L) = \Gamma \oplus
\Psi$ as above, we can write
$$
\sym^2H^0(T,L) \; = \; \sym^2 \Gamma \oplus (\Gamma \otimes \Psi) 
\oplus \sym^2\Psi
$$
with the action of $\tau$ on $\sym^2H^0(T,L)$ having (+1)-eigenspace
$\sym^2 \Gamma\oplus\sym^2\Psi$ of dimension 12 and $(-1)$-eigenspace
$\Gamma \otimes \Psi$ of dimension 9. It follows that the kernel of
the restriction map
$$
\sym^2H^0(T,L) \; \longrightarrow \; H^0(T,L^2)
$$
-- that is, the vector space of quadrics in $\P^5$ vanishing on $T
\subset \P^5$ -- must be contained in the direct sum $\sym^2
\Gamma\oplus\sym^2\Psi$. In other words, we can choose homogeneous
coordinates
$$
[Z,W] \; = \; [Z_0,Z_1,Z_2,W_0,W_1,W_2]
$$
on $\P^5$ so that the action of $\tau$ is given by
$$
\tau \; : \; [Z_0,Z_1,Z_2,W_0,W_1,W_2] \; \mapsto \; 
[Z_0,Z_1,Z_2,-W_0,-W_1,-W_2]
$$
and the defining equations of the double cover $T$ of a general
quartic Enriques surface $S$ may be written in the form
$$
Q_\alpha(Z,W) \; = \; Q'_\alpha(Z) + Q''_\alpha(W).
$$

\

We are now prepared to write down the families of K3 and Enriques
surfaces we'll be studying in the sequel. To start with, let $\Gamma$
and $\Psi$ be 3-dimensional vector spaces, denote by $\P^5$ the
projective space of 1-dimensional subspaces of $\Gamma\oplus\Psi$ And
let
$$
\P^{11} \; = \;\P(\sym^2 \Gamma^\vee\oplus\sym^2\Psi^\vee)
$$
be the projective space of 1-dimensional subspaces of the
(12-dimensional) vector space of quadrics on $\P^5$ of the form above.
Finally, we let $[Z,W] = [Z_0,Z_1,Z_2,W_0,W_1,W_2]$ be homogeneous
coordinates on $\P^5$ with $\Gamma$ the zero locus of $W_0$, $W_1$ and
$W_2$, and $\Psi$ likewise the zero locus of $Z_0$, $Z_1$ and $Z_2$;
and we let $\tau$ be the involution $[Z_0,Z_1,Z_2,W_0,W_1,W_2]
\mapsto[Z_0,Z_1,Z_2,-W_0,-W_1,-W_2]$ of $\P^5$.

\begin{defn}\label{princK3}
By the \textbf{principal family of K3 surfaces} we will mean the
family $\pi : Y \to B$ with $B = \P^{11} \times \P^{11} \times
\P^{11}$ and $Y \subset B \times \P^5$ the subvariety given by
$$
Y \; = \; \left\{\;(Q_1,Q_2,Q_3, [Z,W]) : Q_\alpha(Z,W)=0 \quad 
\forall \alpha=1,2,3 \; \right\},
$$
with $\pi : Y \to B$ the projection on the first factor.
\end{defn}

Note that the action of $\tau$ on the second factor of $B \times \P^5$
carries $Y$ into itself, so that we can make the second

\begin{defn}\label{princenriq}
By the \textit{principal family of Enriques surfaces} we will mean the
family $\pi : X \to B$ with $B$ again as above and $X$ the quotient of
the variety $Y$ above by the involution $\tau$ of $\P^5$.
\end{defn}

It may be a misnomer to call these families of K3 and Enriques
surfaces, since they are only generically that: there are degenerate
fibers, and even fibers of dimension greater than 2.  But it's
convenient to use the term, and we hope the reader will forgive this.

\subsection{Proof of Corollary~\ref{enriqcor}}

In order to apply Theorem~\ref{mainth} to the principal family of
Enriques surfaces and deduce Corollary~\ref{enriqcor}, we simply have
to show that $X \to B$ admits no pseudosections. We'll do this by
analyzing the corresponding family $Y \to B$ of K3 surfaces, since
their equations are in simpler form. We start with the straightforward

\begin{lm}
Let $Y \to B$ be the principal family of K3 surfaces of
Definition~\ref{princK3}. The total space $Y$ is smooth, and its Chow
ring is generated by restrictions of pullbacks of hyperplane classes
under the inclusion
$$
Y \; \hookrightarrow \; \P^{11} \times\P^{11} \times\P^{11} \times \P^5.
$$
\end{lm}

\bpf

To start, introduce the variety
$$
W \; = \;  \{ (Q,p) : p \in Q \} \; \subset \; \P^{11} \times \P^5.
$$
Via the projection $\eta : W \to \P^5$ on the second factor, $W$ is a
$\P^{10}$-bundle over $\P^5$; it's therefore smooth, and its Chow ring
is generated over the Chow ring of $\P^5$ by any class whose
restriction to the fibers of $\eta$ is the hyperplane class on
$\P^{10}$---for example, the restriction of the pullback of the
hyperplane class from $\P^{11}$, via the inclusion $W \hookrightarrow
\P^{11} \times \P^5$. Since the total space $Y$ of our principal
family of K3 surfaces is (via projection to $\P^5$) simply the triple
fiber product
$$
Y \; = \;  W \times_{\P^5} W \times_{\P^5} W
$$
the Lemma follows.

\epf

As an immediate corollary of this Lemma, we have the following
description of cycles $Z \subset X$ of relative dimension 0 over $B$:

\begin{prop}\label{lef}
Let $X \to B$ be the principal family of Enriques surfaces as in
Definition~\ref{princenriq}.  If $Z \subset X$ is any cycle of
codimension 2, the degree of the projection $\pi|_Z : Z \to B$ is
divisible by four.
\end{prop}

\bpf 

Let $\eta: Y \to X$ be the quotient map.  Let $T$ be the class of a
general fiber of $Y$ over $B$.  By the preceding Lemma, the class of
any cycle in $Y$ is a polynomial (with integer coefficients) in the
restrictions to $Y$ of the pullbacks of the hyperplane classes to
$\P^{11} \times\P^{11} \times\P^{11} \times \P^5$. But the first three
of these classes restrict to 0 on a general fiber, so the class of
$\eta^{-1} Z\cdot T$ must be a multiple of the restriction to $T$ of
the hyperplane class on $\P^5$ This has degree divisible by eight.  As
$\eta$ has degree two, the Proposition follows.  

\epf \

As an immediate consequence of Proposition~\ref{lef}, we see that
\textit{the principal family $X \to B$ of Enriques surfaces has no
rational sections}: the image of such a section would give a
codimension 2 cycle of $X$ with degree one over $B$.  

\

In order to show that $X \to B$ admits no pseudosections, it remains
to prove that $X$ cannot contain a subvariety $Z \subset X$ whose
general fiber over $B$ is an irreducible rational curve. To do this,
suppose that $Z$ is such a subvariety.  Let $\tilde Z$ be a resolution
of singularities of $Z$.  We then have a commutative diagram
$$\xymatrix{
\tilde Z \ar[r]^f \ar[rd]^\mu & X\ar[d]^\pi\\
& B}$$

\

Consider the class $f_*(c_1(\omega_\mu))$ in $A^2(X)$.  Since the
general fiber of $\tilde Z$ over $B$ is a smooth rational curve, this
class has degree -2 when restricted to a general fiber of $\pi$.  This
contradicts the fact that all elements of $A^2(X)$ have degree over
$B$ divisible by four.

\

We have thus established the

\begin{lm}
The principal family $X \to B$ of Enriques surfaces admits no
pseudosections.
\end{lm}

\ni Applying Theorem~\ref{mainth} we may deduce
Corollary~\ref{enriqcor}.

\section{Application: torsors for Abelian varieties}\label{Torsors}

It follows from Theorem~\ref{mainth} that any family $\pi:X\rightarrow
B$ of smooth, connected, projective curves of positive genus over some
smooth variety $B$ has a section over $B$ if and only if the
restriction of this family over every curve $C\subset B$ has a
section: since the fibers contain no rational curves, every
pseudosection is a rational section, and every rational section is
everywhere defined.  Similarly, we have the following corollary:

\begin{cor}\label{torsorcor}
Let $B$ be a smooth variety, let $A \rightarrow B$ be an Abelian
scheme over $B$ (i.e. a family of Abelian varieties over $B$), and let
$\pi:T \rightarrow B$ be a torsor for $A\rightarrow B$.  Then $\pi$ is
a trivial torsor if and only if for every curve $C\subset B$, the
restriction $T_C \to C$ is a trivial torsor for $A_C \to C$.
\end{cor}

Since torsors for an Abelian scheme are classified by \'etale
cohomology with coefficients in the group scheme, we can rephrase
Corollary~\ref{torsorcor} by saying that the map
$$ H^1_{\text{\'et}}(B,A) \rightarrow \prod_{C\subset B}
H^1_{\text{\'et}}(C,A_C)$$ 
is injective.  Note that the Lefschetz hyperplane theorem for Picard
groups tells us that the same is true if we replace the Abelian scheme
$A\rightarrow B$ by the commutative group scheme $\mathbb{G}_m \times
B \rightarrow B$.  It is an interesting question, for which other
(possibly noncommutative) group schemes over $B$ does this hold.

\section{Further questions and conjectures}

In this section we will consider some questions and conjectures raised
by Theorem~\ref{mainth}

\subsection{Arithmetic question}

For arithmetic questions related to rationally connected varieties, we
refer the reader to \cite{Abel}.  Let us just mention the following
question which is the analogue of a corollary in the last section.
Let $K$ be a number field, let $B$ be a smooth scheme defined over
$K$, and let $\pi:X \to B$ be a proper, smooth morphism of schemes
whose geometric fibers are connected curves of positive genus.
Suppose that for every number field extension $L/K$, the induced
mapping on rational points $\pi:X(L) \to S(L)$ is surjective -- we
refer to this property by saying $\pi$ is \emph{arithmetically
surjective}.  Does it then follow that $\pi:X \to S$ has a section?
We may also ask the same question when the geometric fibers of $\pi$
are Abelian varieties.

\subsection{Possible extensions}

We would like to take a moment here to discuss possible extensions of
Theorem~\ref{mainth}. To begin with, we interpreted the theorem as
stating that a family of varieties $\pi : X \to B$ such that every
one-parameter subfamily has a section has this property ``by virtue
of" the fact that $X$ contains a family of rationally connected
varieties. But the statement of the theorem asserts only the existence
of a pseudosection in $\pi : X \to B$; it doesn't assert any direct
connection between the sections of $X_C \to C$ over very general
curves $C$ and the pseudosection.  Accordingly, we could ask:

\begin{ques} 

Does there exist a family $\cH'_d $ of curves on $B$, whose general
member is smooth and irreducible, with the property that for any
proper morphism $\pi : X \to B$ of relative dimension $d$, for very
general $[C] \in \cH'_d$ every section of the pullback
$$
 \pi_C : X_C = X \times_B C \to C
$$
lies in a pseudosection of $\pi : X \to B$?  

\end{ques}

One special case of this question is when $\pi:A \rightarrow B$ is an
Abelian scheme over a smooth variety $B$.  In this case we are asking
whether we can find a family of curves $C\subset B$ such that for a
very general member of this family, the map
$$ H^0_{\text{\'et}}(B,A) --> H^0_{\text{\'et}}(C,A_C) $$
is surjective.  So, in this case, the question above is an
$H^0$-analogue of the $H^1$-interpretation of Theorem~\ref{mainth}.

\subsection{Dependence on $d$}

A simpler question is whether we can eliminate the dependence of the
family $\cH_d$ of curves on $d$. The answer to this seems to be ``no,"
and while we cannot prove the impossibility it seems worthwhile to
describe here the examples that lead us to this conclusion, since it
may shed some light on how fast the size of the curves in $\cH_d$ have
to grow with $d$.

\

Briefly, for any $e$ we will write down families of hypersurfaces in
$\P^n$ parametrized by $B = \P^2$ with the property that their
restriction to any curve $C \subset \P^2$ of degree $e$ or less has a
section, but which we believe to admit no pseudosections at all. To do
this, let $m$ and $n$ be any integers; let $\P^n$ be projective
$n$-space with homogeneous coordinates $[X_0,\dots,X_n]$ and let
$\P^N$ be the projective space parametrizing hypersurfaces of degree
$m$ in $\P^n$. Let $\cX \subset \P^N \times \P^n \to \P^N$ be the
universal hypersurface of degree $m$ in $\P^n$; that is, the zero
locus of the polynomial
 $$
 F(a,X) \; = \; \sum a_IX^I
 $$
which is linear in the coordinates $a_I$ on $\P^N$ and of degree $m$
in the $X_i$. Finally, let $\P^2 \hookrightarrow \P^N$ be a general
map of the form
 $$
 \P^2 \; \xrightarrow{\nu_e} \;  \P^{\binom{e+2}{2}-1} \; \to \; \P^N
 $$
where $\nu_e$ is a Veronese map of degree $e$ and the second map is a
general linear inclusion; and let
 $$
 \pi \; : \; X = \P^2 \times_{\P^N}  \cX \; \longrightarrow \; \P^2
 $$
be the pullback of the universal hypersurface to $\P^2$ via this inclusion.

Now assume that
 $$
 \binom{e+2}{2} \; = \; n+1
 $$
and that $m$ is large.  Consider the following two assertions:

\

\ni $\bullet$ \thinspace The restriction of the family $\pi : X \to
\P^2$ to any curve $C \subset \P^2$ of degree $e$ or less has a
section; but

\

\ni $\bullet$ \thinspace The family $\pi : X \to \P^2$ itself has no
pseudosection.

\

The first of these assertions is straightforward to prove: under the
inclusion $\P^2 \hookrightarrow \P^N$, the span of a curve $C \subset
\P^2$ of degree $e$ or less has dimension $\binom{e+2}{2} - 1 = n$ or
less. Thus the hypersurfaces appearing as fibers of the restriction
$X_C \to C$ of the family $\pi : X \to \P^2$ to $C$ are all linear
combinations of $n$ hypersurfaces $G_1,\dots,G_n \subset \P^n$, and
any point of intersection of these hypersurfaces gives a section of
$X_C \to C$.

\

As for the second assertion, we cannot prove it but we give a
``plausibility argument'' which suggests it is true.  To begin with, a
general fiber of $\pi : X \to \P^2$ is a general hypersurface of
degree $m$ in $\P^n$; by a result of Clemens \cite{clemens}, for $m$
large this will contain no rational curves. Thus to prove the second
assertion we need only show that $\pi : X \to \P^2$ has no rational
sections.

\

Since rational sections over $\P^2$ are tricky to parametrize we will
restrict to a general curve $C \subset \P^2$ of degree $e+1$, and
present evidence that the restriction $X_C \to C$ has no section.  To
do this, we start by counting the dimension of the family of sections
of the product $C \times \P^n$ there are of a given degree $k$---that
is, graphs of maps $C \to \P^n$ of degree $k$---and then estimating
the number of conditions it imposes on such a section to require it
lies on the hypersurface $X_C \subset C \times \P^n$. For the first, a
map $C \to \P^n$ of degree $k$ is given by a line bundle $L$ of degree
$k$ on $C$, together with $n+1$ sections of $L$ up to scalars. The
line bundles of degree $k$ on $C$ are parametrized by the Jacobian of
$C$, which has dimension
 $$
 g \; = \; \binom{e}{2}.
 $$
If $k$ is large, moreover, each such line bundle will have $k-g+1$
global sections, so the dimension of the family of maps $C \to \P^n$
of degree $k$ is
 $$
 g + (n+1)(k-g+1) - 1 \; = \; (n+1)(k+1) - ng - 1.
 $$

\

Now let's count how many conditions it is for the graph of such a map
to lie in $X_C$. This is straightforward: when we pull the polynomial
$F(a,X)$ defining the universal hypersurface back to $C$, the
coefficients pull back to section of $\cO_C(e)$ and the coordinates
$X_i$ to sections of $L$, so that the pullback of $F$ is a section of
the bundle
 $$
 M \; = \; L^{\otimes m} \otimes \cO(e).
 $$
The number of conditions for this section to vanish identically should
thus be 
\begin{align*} h^0(M) \; &= \; \deg(M) - g + 1 \\ &= \; km +
e(e+1) - g + 1 \end{align*} 
and the expected dimension of the family
of sections of $X_C \to C$ of degree $k$ is accordingly
 $$
   (n+1-m)k - (n-1)(g-1) - e(e+1).
 $$
In particular, for $m$ large this is negative, suggesting that there 
should be no sections.


\begin{thebibliography}{[GHMS]}

\bibitem[Ca]{Ca} F. Campana, {\it Connexit\'e rationnelle des 
vari\'et\'es de Fano}, Ann. Sc. E.N.S. {\bf25} (1992), 539-545

\bibitem[C]{clemens} H. Clemens, {\it Curves on generic hypersurfaces},
Ann. Sci. \'Ecole Norm. Sup. (4) {\bf 19} (1986) 629-636.

\bibitem[FP]{FP} W. Fulton, R. Pandharipande, {\it Notes on stable
    maps and quantum cohomology}, in {\em Algebraic geometry -- Santa
    Cruz 1995}, AMS (1995), 45-96.

\bibitem[GHMS]{Abel} T. Graber, J. Harris, B. Mazur, J. Starr, {\it
Arithmetic questions related to rationally connected varieties}, in
preparation.

\bibitem[GHS]{GHS} T. Graber, J. Harris, J. Starr, {\it Families
of rationally connected varieties} to appear in J. Amer. Math. Soc.

\bibitem[GS]{GS} A. Grothendieck and J.-P. Serre, {\em Correspondance
Grothendieck-Serre}, ed. by Pierre Colmez and J.-P. Serre, 
Soci\'et\'e Math. de France 2001.

\bibitem[GM]{gm} M. Goresky and R. MacPherson, {\em Stratified Morse
Theory}, Ergebnisse der Math. 14, Springer-Verlag, Berlin, 1988. 

\bibitem[K]{K} J. Koll\'ar,  {\em Rational Curves on Algebraic 
Varieties}, Ergebnisse der Math.
32, Springer-Verlag, Berlin, 1996.

\bibitem[KMM]{KMM}  J. Koll\'ar, Y. Miyaoka, S. Mori, {\em Rationally 
Connected Varieties}, J.
Alg. Geom. {\bf 1}  (1992) 429-448.

\end{thebibliography}
\end{document}